\documentclass{elsart}
\usepackage{graphicx}
\usepackage{latexsym}
\usepackage{amssymb}
\usepackage{subfigure}

\usepackage{graphicx,epsfig,amssymb,amsmath,amsfonts}
\journal{Communications in Nonlinear Science and Numerical Simulation}

\begin{document}
\begin{frontmatter}
\title{Numerical generation of periodic traveling wave solutions of some nonlinear dispersive wave equations}
\author{J. \'Alvarez}
\address{Department of Applied Mathematics,
University of Valladolid, Paseo del Cauce 59, 47011,
Valladolid, Spain.}
\address{
IMUVA, Institute of Mathematics of University of Valladolid; Spain.
Email: joralv@eii.uva.es}
\author{A. Dur\'an \thanksref{au}}
\address{Department of Applied Mathematics, University of
Valladolid, Paseo de Belen 15, 47011-Valladolid, Spain.}
\address{
IMUVA, Institute of Mathematics of University of Valladolid; Spain.
Email:
angel@mac.uva.es }
\thanks[au]{Corresponding author}

\begin{abstract}
Proposed in this paper is a numerical procedure to generate periodic traveling wave solutions of some nonlinear dispersive wave equations. The method is based on a suitable modification of a fixed point algorithm of Petviahvili type and solves several drawbacks of some previous algorithms proposed in the literature. The method is illustrated with the numerical generation of periodic traveling waves of fractional KdV type equations and some extended Boussinesq systems.
\end{abstract}
\end{frontmatter}
\section{Introduction}
\label{se1}
Considered here is a technique to generate numerically periodic traveling wave solutions of nonlinear wave equations of the form
\begin{eqnarray}
u_{t}-Lu_{x}+f(u)_{x}=0,\quad x\in\mathbb{R}, t\geq 0.\label{ptw1}
\end{eqnarray}
In (\ref{ptw1}) $u=u(x,t)$, $f:\mathbb{R}\rightarrow\mathbb{R}$ is a smooth function (typically a polynomial) with $f(0)=f^{\prime}(0)=0$ and $L$ is a dispersion operator with Fourier symbol $\alpha: \mathbb{R}\rightarrow\mathbb{R}$,
\begin{eqnarray}
\widehat{Lv}(\xi)=\alpha(\xi)\widehat{v}(\xi),\quad \xi\in\mathbb{R},\label{ptw2}
\end{eqnarray}
being $\widehat{v}(\xi)$ the Fourier transform of $v$ at $\xi\in\mathbb{R}$. Equation (\ref{ptw1}) appears as a first one-dimensional, unidirectional approach for nonlinear dispersive wave propagation, see e.~g. \cite{BenjaminBM1972,AbdelouhabBFS1989} and particular cases are remarkable: Korteweg de Vries-type equation (KdV, $L=-\partial_{xx}, f(u)=u^{p-1}, p\geq 3$), Benjamin-Ono-type equation (BO, $L=-\mathcal{H}\partial_{x}$ with $\mathcal{H}$ the Hilbert transform, $f(u)=u^{p-1}, p\geq 3$) or Benjamin-type equation (BO, $L=-\beta\partial_{xx}-\gamma\mathcal{H}\partial_{x}$ with $\beta,\gamma\geq 0, f(u)=u^{p-1}, p\geq 3$).

The paper is concerned with the numerical generation of periodic traveling wave solutions of the models (\ref{ptw1}). They are solutions of the form
\begin{eqnarray}
u(x,t)=\phi(x-c_{s}t)=\sum_{n=-\infty}^{\infty} \widehat{\phi}_{n}e^{i\frac{n\pi x}{l}(x-c_{s}t)},\label{ptw3}
\end{eqnarray}
whic are periodic of some period $2l, l>0$ and travel with speed $c_{s}>0$. Results on existence and stability of such solutions for particular cases of (\ref{ptw1}) can be seen in the literature (see e.~g. \cite{Angulo2006} and references therein). Recently, \cite{ChenB2013}, Chen and Bona extablished existence of periodic traveling wave solutions of (\ref{ptw1}) under some hypotheses on $f$ and $L$ (see below) as well as stability under perturbations of the same period. This study generalizes previous works on this matter, \cite{Chen2004}.

In \cite{AlvarezD2011} a numerical scheme for the numerical generation of periodic traveling waves of (\ref{ptw1}) was proposed. This was based on a Fourier pseudospectral discretization in space of (\ref{ptw1}) along with time integration with conservation properties of invariant quantities of the periodic value problem associated to (\ref{ptw1}). The periodic traveling wave profile was approximated by using the iterative Petviashvili method, \cite{petviashvili,AlvarezD2014a}. The resulting scheme ensured an accurate computation of two characteristic elements of the waves, the initial profile and the speed. The present work attempts to complete \cite{AlvarezD2011} in two main points, not totally resolved there in the authors' opinion. The first one is that the generation of the profile of the periodic traveling wave solution may involve some nonzero constants which make the Petviashvili method, at least initially,  unable to be used. The second one, somehow related, is that this method is not well suited to nonlinearities $f$ in (\ref{ptw1}) with more than one homogeneous terms of different degrees (like a polynomial for example).

In order to overcome these drawbacks, the present work introduces an improved algorithm from the one proposed in \cite{AlvarezD2011}. The strategy is based on a change of variable suggested in \cite{ChenB2013} (see also \cite{ChenCN2007}) to treat the constants (to solve the first problem above described) and the use of extended versions of the Petviashvili method introduced in \cite{AlvarezD2014b} that allows the presence of more general nonlinearities in the model (to solve the second problem). The structure is as follows. In Section \ref{se2} some hypotheses on (\ref{ptw1}) (considered in \cite{ChenB2013}) are assumed and the proposed algorithm is decribed. The efficiency of the method is checked in Section \ref{se3} by considering the periodic traveling wave computation for two models. The first one is the generalized fractional KdV equation (of the form (\ref{ptw1}) with $\alpha(\xi)=|\xi|^{\mu}, 0<\mu\leq 2, f(u)=u^{p-1}/(p-1)$), \cite{AbdelouhabBFS1989,Johnson2013}. The second model is the e-Boussinesq system, \cite{NguyenD2008}. This is not of the form (\ref{ptw1}) but its presence here is justified as an attempt to extend the application of the method is based to more general situations, in particular to some nonlinear dispersive systems. In both cases and to our knowledge, many of the profiles are generated for the first time; in particular, for the case of the e-Boussinesq system, we do not even know if some theoretical results of existence of periodic traveling wave solutions for the model are available (although they may be expected from the existence of solitary wave solutions as their large wavelength limit, \cite{Bona1981}). Section \ref{se4} contains some concluding remarks.
\section{Preliminaries and numerical method}
\label{se2}
We make the following assumptions about (\ref{ptw1}):
\begin{itemize}
\item[(H1)] The nonlinear term $f$ is a polynomial of the form
\begin{eqnarray*}
f(z)=3\gamma_{3}z^{2}+\cdots +p\gamma_{p}z^{p-1},\label{ptw4}
\end{eqnarray*}
where $\gamma_{j}\geq 0, j=3,\ldots p-1, \gamma_{p}>0, p\geq 3$.
\item[(H2)] $\alpha$ is real, even, continuous with $\alpha(0)=0$. If $p_{0}=\min\{j=3,\ldots,p/ \gamma_{j}>0\}$ then there exists $\widetilde{s}\geq \displaystyle\frac{p_{0}-2}{4}$ such that
$$\lim_{\xi\rightarrow 0}\frac{\alpha(\xi)}{|\xi|^{2\widetilde{s}}}=0.$$
\item[(H3)] There exists ${s}\geq \displaystyle\frac{p-2}{4}$ such that
$$0<\lim_{\xi\rightarrow \infty}inf \frac{\alpha(\xi)}{|\xi|^{2{s}}}\leq 
\lim_{\xi\rightarrow \infty} sup \frac{\alpha(\xi)}{|\xi|^{2{s}}}<\infty.$$
\end{itemize}
These are some of the hypotheses considered in \cite{ChenB2013} and will be adopted here with the aim of ensuring the existence of periodic traveling wave solutions in the following sense (see \cite{ChenB2013} for details): For a sufficiently large period parameter $l>0$ the profiles $\phi$ are obtained as minimizers of constrained variational problems involving momentum and energy functionals for (\ref{ptw1}). Then for each minimizer $\phi$ one can find a positive speed $c_{s}$ such that (\ref{ptw3}) is an infinitely smooth traveling wave solution of (\ref{ptw1}) which is also stablñe under perturbations of the same period. In \cite{ChenB2013} the authors also study the large wavelength limit $l\rightarrow\infty$ of the periodic traveling waves determining under some additional hypotheses the convergence to solitary wave solutions.

A comment in \cite{ChenB2013} is the starting point to introduce our numerical strategy to generate numerically periodic traveling wave profiles. When $u$ in (\ref{ptw3}) is a solution of (\ref{ptw1}) then the profile $\phi$ must satisfy
\begin{eqnarray}
-c_{s}\phi-L\phi+f(\phi)=A,\label{ptw5}
\end{eqnarray}
where $A$ is a (nonzero in general) constant. Then our method is based on the following steps:
\begin{itemize}
\item[(S1)] Finding (approximate or exact) real constants solutions $\phi=C$ of (\ref{ptw5}) that is
\begin{eqnarray}
-c_{s}C+f(C)=A.\label{ptw6}
\end{eqnarray}
(Note that $LC=0$ since $\alpha(0)=0$ is assumed in (H2).) This means that $C$ must be a root of the polynomial
$$P(z)=-A-c_{s}z+3\gamma_{3}z^{2}+\cdots +p\gamma_{p}z^{p-1}.$$
\item[(S2)] The change of variables $\varphi=\phi-C$, when $\phi$ is a solution of (\ref{ptw5}) along with (\ref{ptw6}) lead to the following equation for $\varphi$:
\begin{eqnarray}
\underbrace{(c_{s}+L-f^{\prime}(C))}_{\widetilde{L}(C)}\varphi-\underbrace{\sum_{j=2}^{p-1}\frac{f^{j)}(C)}{j!}\varphi^{j}}_{N(\varphi)}=0,\label{ptw7}
\end{eqnarray}
where
\begin{eqnarray*}
f^{\prime}(C)&=&\sum_{l=3}^{p}l(l-1)\gamma_{l}C^{l-2}\\
f^{j)}(C)&=&\sum_{l=j+1}^{p}l(l-1)\cdots (l-j)\gamma_{l}C^{l-j-1},\quad j\geq 2.
\end{eqnarray*}
\item[(S3)] Equation (\ref{ptw7}) is iteratively solved with a extended version of the Petviashvili method, \cite{AlvarezD2014b}. Writing
$$
N(\varphi)=\sum_{j=2}^{p-1}N_{j}(\varphi),\quad N_{j}(\varphi)=\frac{f^{j)}(C)}{j!}\varphi^{j}, j=2,\ldots,p-1,
$$ and from $\varphi^{[0]}\neq 0$  and for $\nu=0,1,\ldots$ the following iteration is considered:
\begin{eqnarray}
  \widetilde{L}(C)\varphi^{[\nu+1]}&=&\sum_{j=2}^{p-1}s_{j}(\varphi^{[\nu]})N_{j}(\varphi^{[\nu]}), \nu=0,1,\ldots,\label{ptw8}\\
s_{j}(u)&=&s(u)^{\alpha_{j}},\quad s(u)=\left(\frac{\langle \widetilde{L}(C)u,u\rangle}{\langle N(u),u\rangle}\right),\label{ptw9}\\
&&\quad \alpha_{j}=\frac{{j}}{{j}-1},\quad j=2,\ldots,p-1,\nonumber
\end{eqnarray}
The factors $s_{j}$ in (\ref{ptw9}) are obtained by using the same quotient $s(u)$ but changing the exponent $\alpha_{j}$. Formulas (\ref{ptw8}), (\ref{ptw9}) correspond to an extension of the classical Petviashvili scheme for nonlinearities with several homogeneous terms but different degree of homogeneity. This belongs to a more general family, described in \cite{AlvarezD2014b}. The use of (\ref{ptw8}), (\ref{ptw9}) instead of other alternatives presented in the literature (of Newton-type, for example) is justified by the attempt to simplify the implementation of the iterative techniques for the profiles with an easily implementable fixed point type algorithm, \cite{yang2}. Note that in the case of convergence, $s(\varphi^{[\nu]})\rightarrow 1$.
\end{itemize}
The operator (\ref{ptw2}) taking part of $\widetilde{L}(C)$ in (\ref{ptw7}) suggests to consider the implementation of (\ref{ptw8}), (\ref{ptw9}) in Fourier space for the Fourier coefficients $\widehat{\varphi^{[\nu]}}_{n}, n\in \mathbb{Z}$ of the iterations $\varphi^{[\nu]}, \nu=0,1,\ldots$ Then for fixed $c_{s}$ and $l>0$ the problem (\ref{ptw5}) on $(-l,l)$ with periodic boundary conditions is considered and (\ref{ptw8}), (\ref{ptw9}) take the form
\begin{eqnarray}
\widehat{\varphi^{[\nu+1]}}_{n}&=&\frac{1}{c_{s}+\alpha\left(\frac{n\pi}{l}\right)-f^{\prime}(C)}
\sum_{j=2}^{p-1}s_{j}(\varphi^{[\nu]})\widehat{\left(N_{j}(\varphi^{[\nu]})\right)_{n}}\label{ptw8b}\\
 &&n\in \mathbb{Z},\quad \nu=0,1,\ldots,\nonumber \\
s_{j}(u)&=&\left(\frac{\sum_{n\in\mathbb{Z}}
(c_{s}+\alpha\left(\frac{n\pi}{l}\right)-f^{\prime}(C))|\widehat{u}_{n}|^{2}}{\sum_{n\in\mathbb{Z}}
\widehat{\left(N(u)\right)}_{n}\overline{\widehat{u}}_{n}}\right)^{\alpha_{j}},\label{ptw9b}
\\
&&\quad j=2,\ldots,p-1,\nonumber
\end{eqnarray}
%

Finally, we note that the existence results in \cite{ChenB2013} are mentioned here in connection with the hypotheses (H1)-(H3). However, the strategy for the numerical generation uses a different argument, in the sense that it previously fixes a speed $c_{s}$ whose associated profile is attempted to be computed. In this sense, the point of view of the resolution is closer to other existence results presented in the literature, \cite{AnguloBS2006}.
\section{Numerical examples}
\label{se3}
This section is devoted to illustrate the behaviour of the algorithm to generate numerically periodic traveling waves. This will be shown by considering two models: a generalized version of the fractional KdV equation and the extended Boussinesq system.

It may be worth making first some brief comments on the implementation of (\ref{ptw8}), (\ref{ptw9}). As mentioned above, due to the nonlocal operator (\ref{ptw2}), a natural way to approximate (\ref{ptw7}) consists of a discretization on an interval with a Fourier pseudospectral method as in \cite{AlvarezD2011}, see \cite{Boyd2000}. This allows to implement the iteration by using the discrete Fourier version of (\ref{ptw8b}), (\ref{ptw9b}) and FFT techniques. Secondly, in some of the experiments, the resulting scheme has been complemented with vector extrapolation methods to accelerate the convergence. The literature on these methods is abundant (see e.~g. \cite{BrezinskiR1991,Sidi2003} and references therein) and their use in traveling wave simulations is analyzed in  \cite{AlvarezD2015b}.
\subsection{Generalized fractional KdV equation}
The first example treated here concerns the generation of periodic traveling waves of fractional KdV models of the form
\begin{eqnarray}
u_{t}-\Lambda^{\mu}u_{x}+(f(u))_{x}=0,\label{ptw10}
\end{eqnarray}
which corresponds to (\ref{ptw1}) with $L=\Lambda^{\mu}$, being 
$$\widehat{\Lambda v}(\xi)=|\xi|\widehat{v}(\xi),$$ $0<\mu\leq 2$ and $f(u)=\displaystyle\frac{u^{p-1}}{p-1}, p\geq 3$. Existence and stability of periodic traveling wave solutions are analyzed in \cite{Johnson2013} by using variational arguments. Specifically, existence is obtained for $0<p<p_{max}$ where
$$
p_{max}=\left\{\begin{matrix}\frac{2}{1-\mu}&\mbox{if}\quad \mu<1\\+\infty&\mbox{if} \quad \mu\geq 1\end{matrix}\right.
$$
As far as stability is concerned, for $\mu>1/2$ and
$$p^{*}(\mu)=\frac{(3+\mu)2^{\mu}-2\mu+(\mu-1)2^{\mu+1}}{2+(\mu-1)2^{\mu}},$$
then the periodic traveling wave is 
\begin{itemize}
\item Spectrally stable if $\mu>1$ and $1\leq p<p^{*}(\mu)$.
\item Spectrally unstable if $1/2<\mu<1$ or $\mu>1$ when $p>p^{*}(\mu)$.
\end{itemize}
Furthermore, if 
$$p^{*}_{max}=\max_{\mu\geq 1} p^{*}(\mu),$$
then for $p\in (1,p^{*}_{max})$ there are $\mu_{-}(p)<\mu_{+}(p)$ such that small amplitude periodic traveling waves are spectrally stable for $\mu\in (\mu_{-}(p),\mu_{+}(p))$ and spectrally unstable for $\alpha\notin (\mu_{-}(p),\mu_{+}(p))$.

\begin{figure}[htbp]
\centering 
\subfigure[]{
\includegraphics[width=6.6cm]{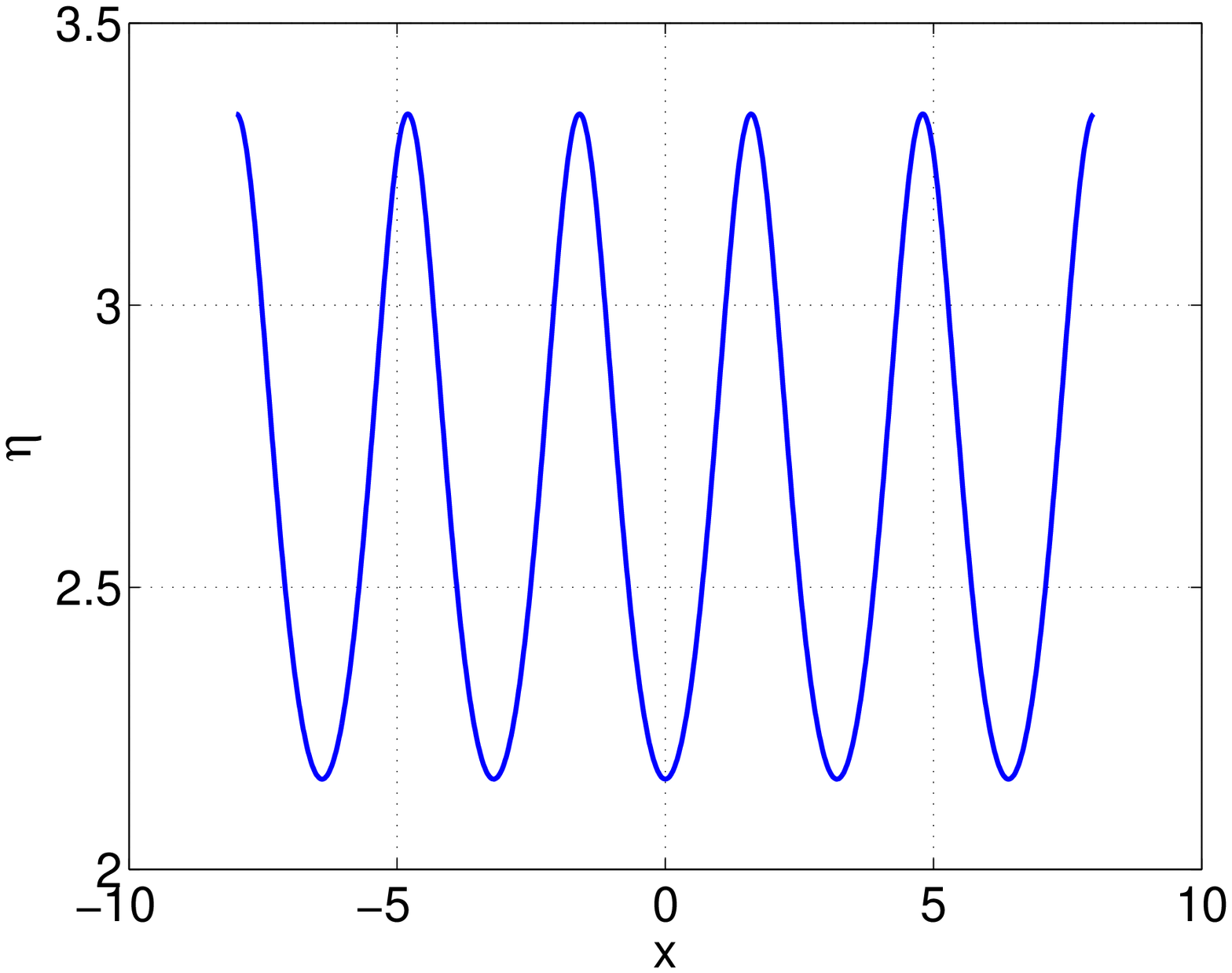} }
\subfigure[]{
\includegraphics[width=6.6cm]{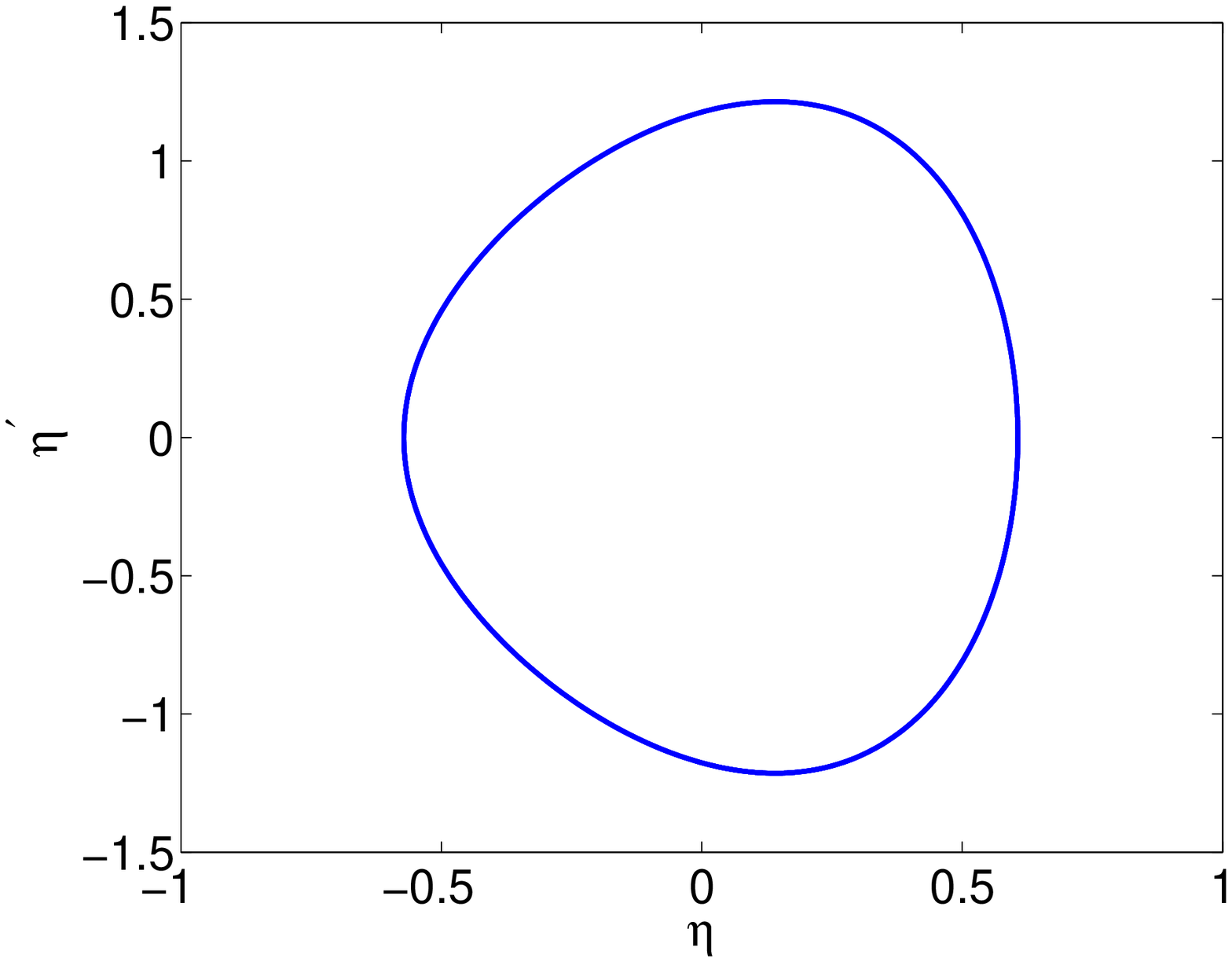} }
\subfigure[]{
\includegraphics[width=6.6cm]{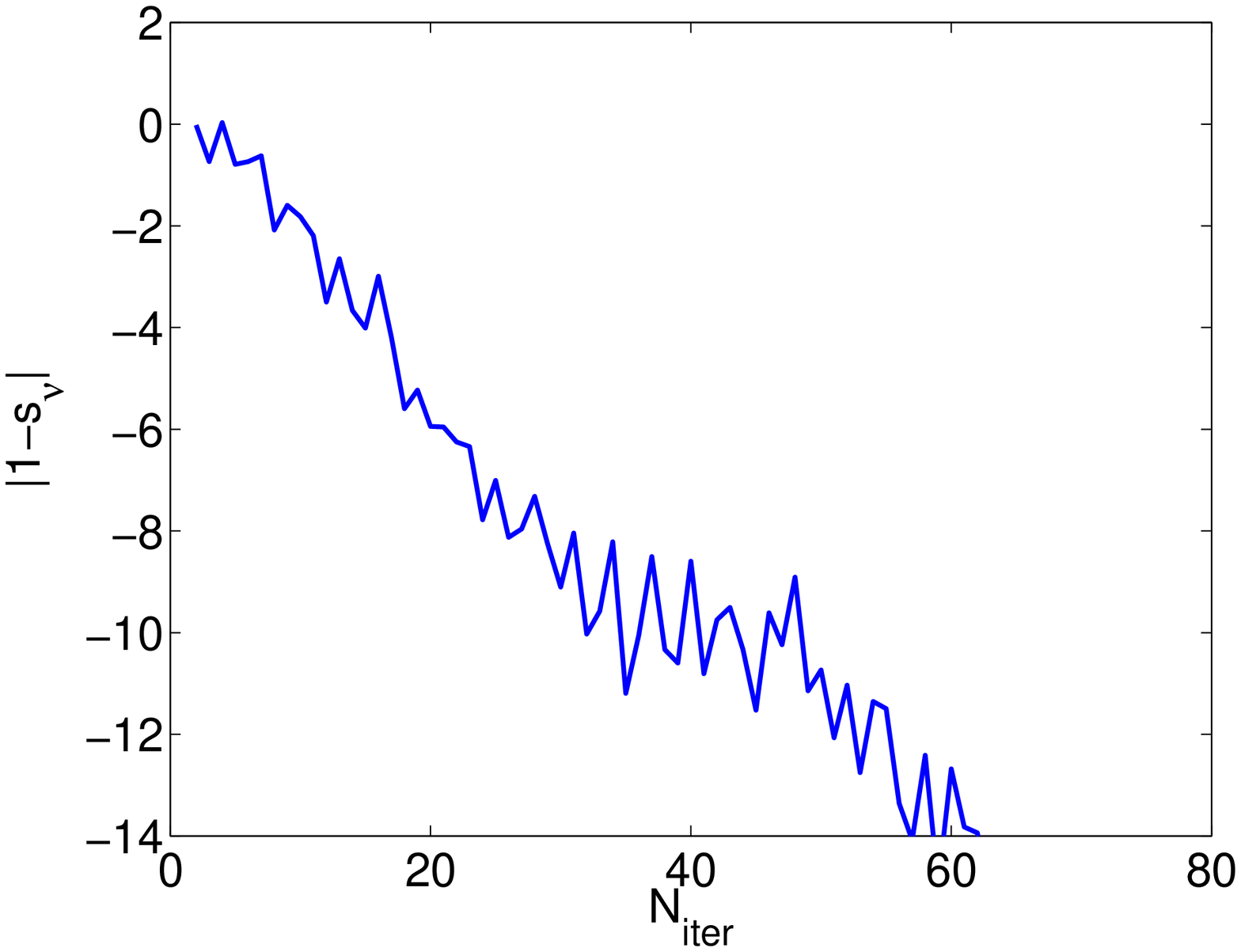} }
\subfigure[]{
\includegraphics[width=6.6cm]{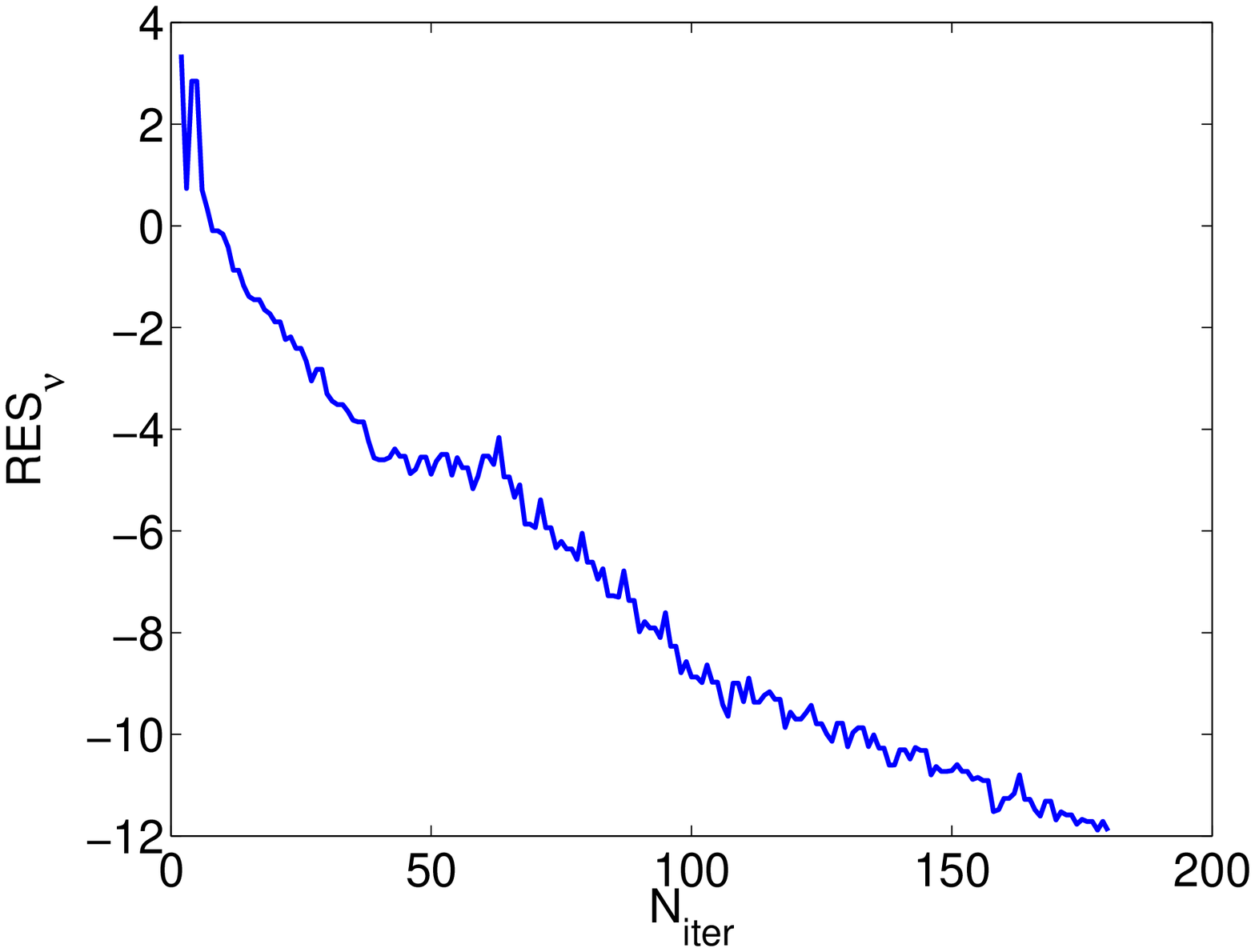} }
\caption{Periodic traveling wave generation of  (\ref{ptw10}) with $A=c_{s}=1$ and $p=3,\alpha=0.8$. (a) Approximate profiles; (b) Phase portrait; c) Discrepancy (\ref{fsec33})  vs number of iterations; (d)
Residual error (\ref{fsec32}) vs number of iterations. (Semi-logarithm scale in both cases.)  
} \label{figptw1}
\end{figure}

\begin{figure}[htbp]
\centering \subfigure[]{
\includegraphics[width=6.6cm]{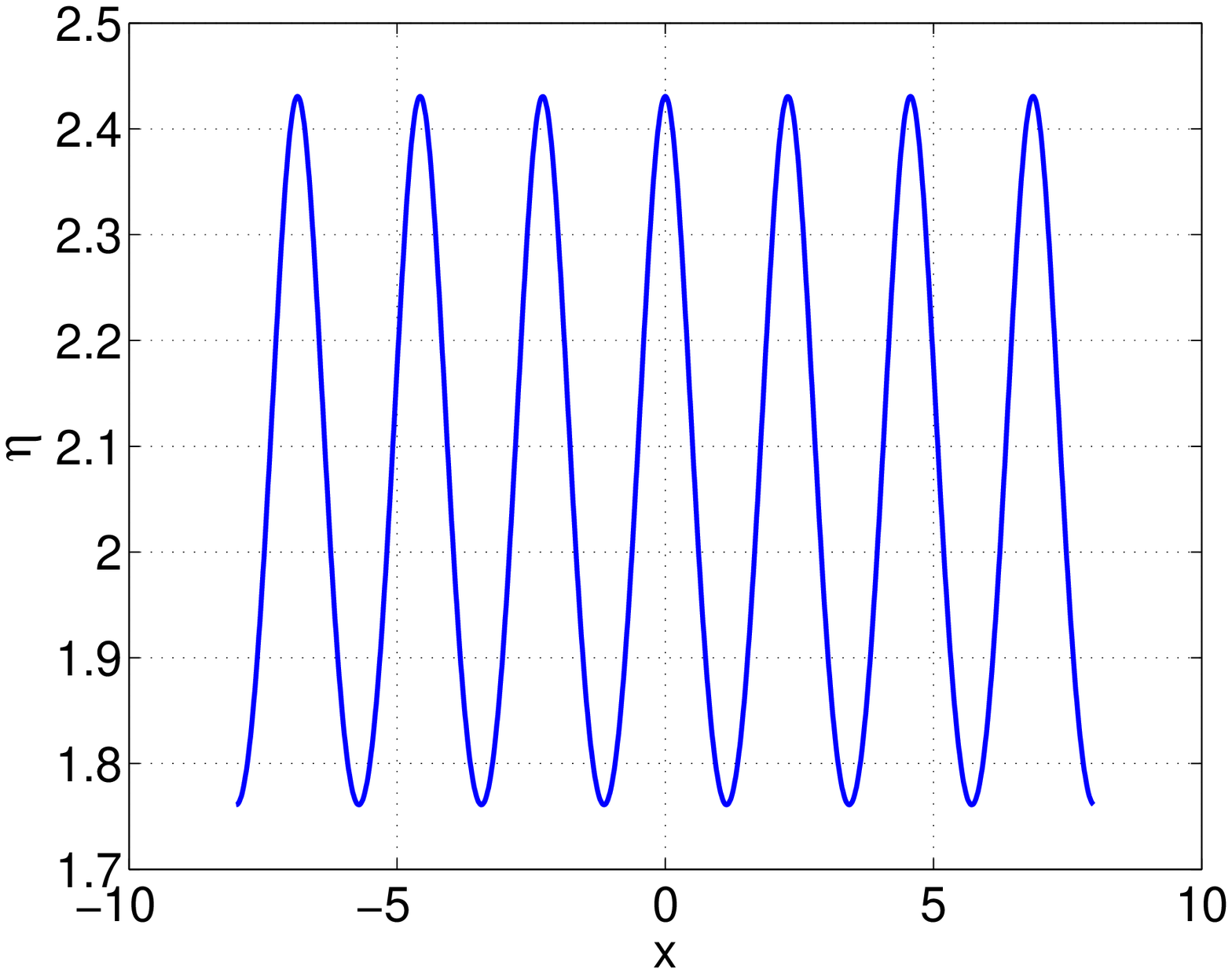} }
\subfigure[]{
\includegraphics[width=6.6cm]{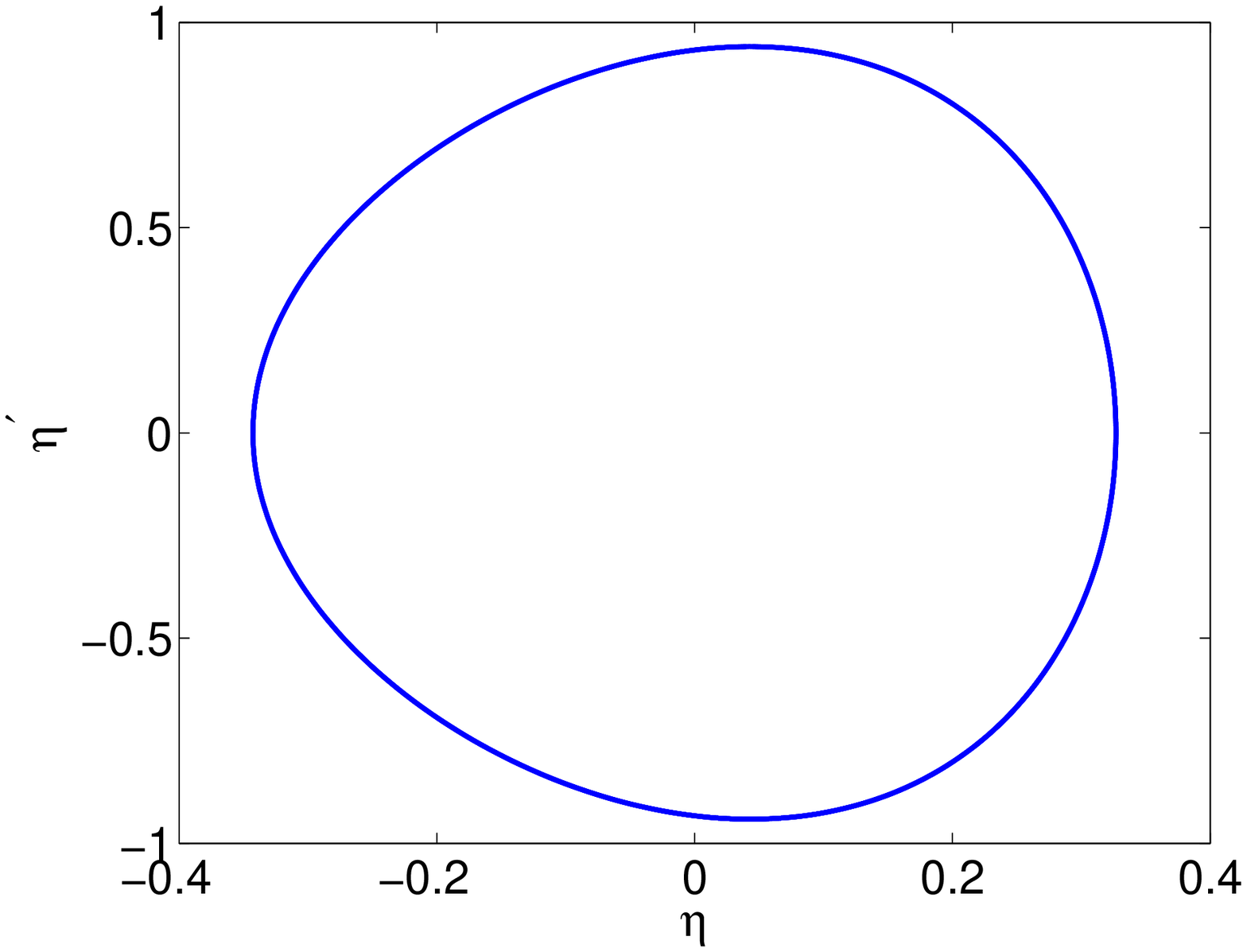} }
 \subfigure[]{
\includegraphics[width=6.6cm]{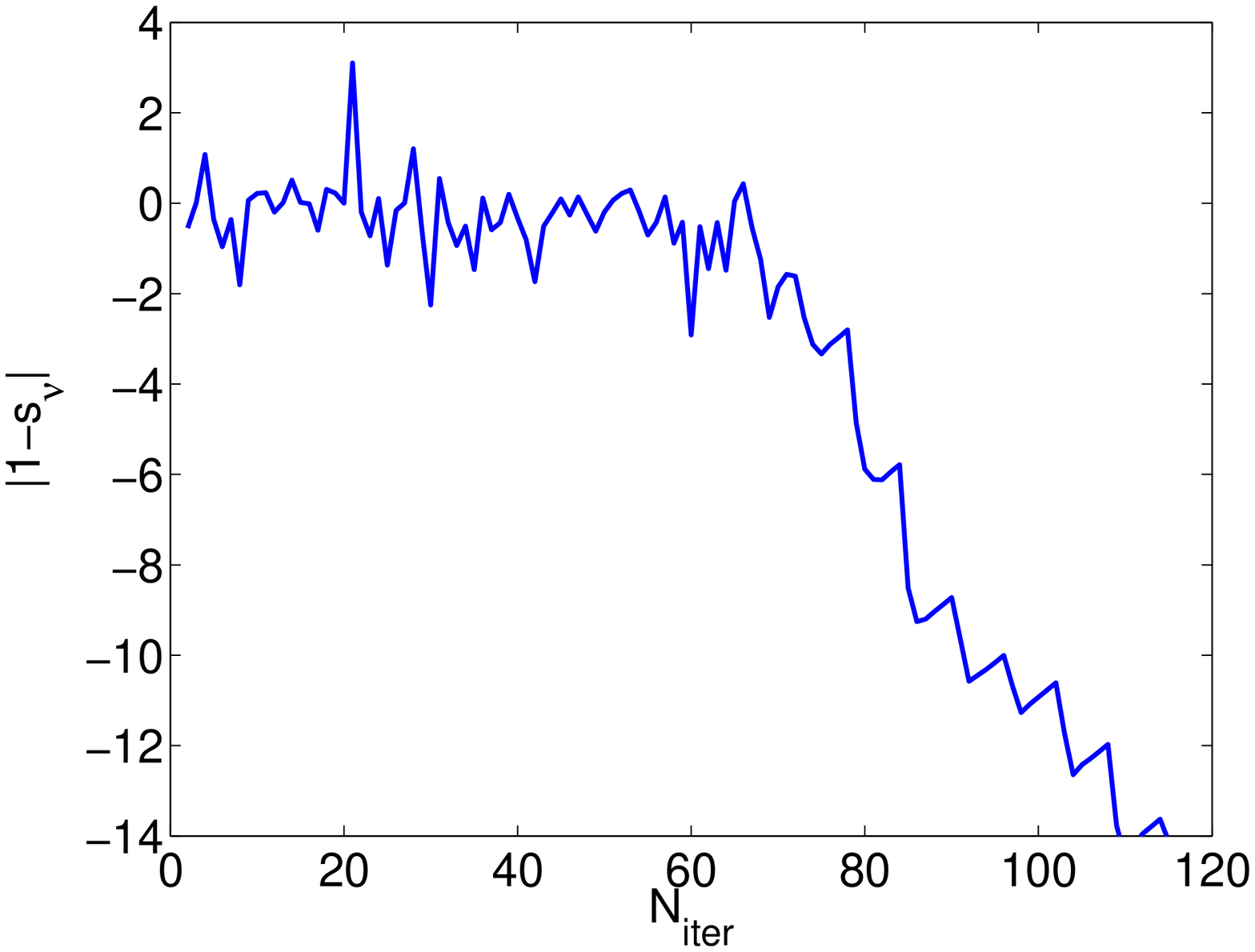} }
\subfigure[]{
\includegraphics[width=6.6cm]{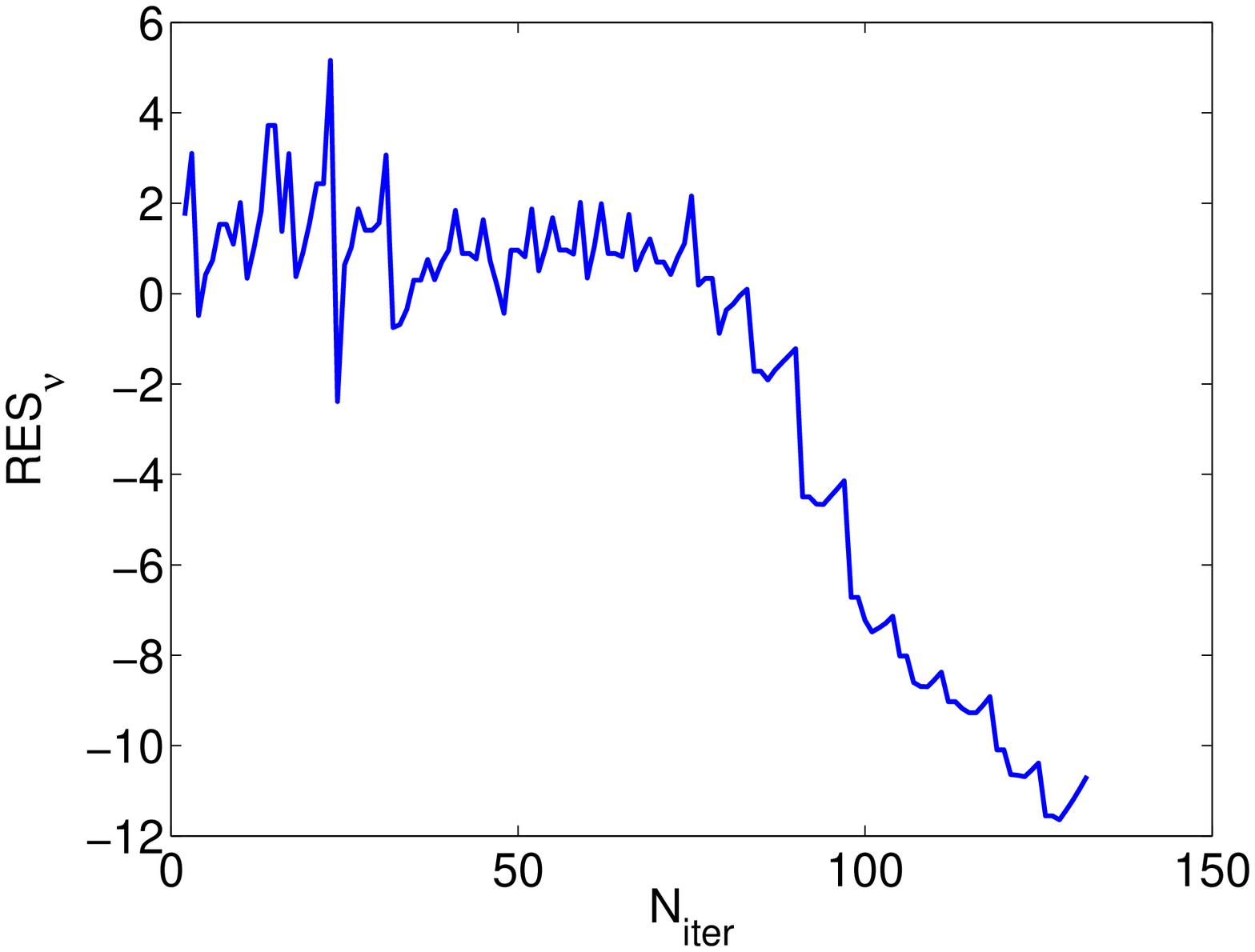} }
\caption{Periodic traveling wave generation of  (\ref{ptw10}) with $A=c_{s}=1$ and $p=4,\alpha=1.2$. (a) Approximate profiles; (b) Phase portrait; c) Discrepancy (\ref{fsec33})  vs number of iterations; (d)
Residual error (\ref{fsec32}) vs number of iterations. (Semi-logarithm scale in both cases.) } \label{figptw2}
\end{figure}

In this case step [S1] of the algorithm requires the computation of constant solutions $u=C$ satisfying
\begin{eqnarray}
\frac{C^{p-1}}{p-1}-c_{s}C-A=0,\quad A\in\mathbb{R}.\label{ptw11}
\end{eqnarray}
Existence of real solutions of (\ref{ptw11}) requires some relations involving the parameters $p, c_{s}, A$ (but not $\mu$). From the constant $C$ of (\ref{ptw11}) the equation (\ref{ptw7}) takes the form 
\begin{eqnarray*}
\underbrace{(c_{s}+L-C^{p-2})}_{\widetilde{L}(C)}\varphi-\underbrace{\frac{1}{p-1}\sum_{j=2}^{p-1}\begin{pmatrix}p-1\\ j\end{pmatrix}C^{p-1-j}\varphi^{j}}_{N(\varphi)}=0,\label{ptw12}
\end{eqnarray*}
The behaviour of the algorithm is illustrated here by fixing $A=c_{s}=1$ and two values for $\alpha$ and $p$: $\alpha=0.8, p=3$ and $\alpha=1.5, p=4$. In both cases, the convergence was accelerated by using the minimal polynomial extrapolation method (MPE), \cite{BrezinskiR1991,Sidi2003}. The approximate profile along with the corresponding phase portrait, are shown in Figures \ref{figptw1}(a),(b) and \ref{figptw2}(a),(b), respectively.
Figures \ref{figptw1}(c),(d) and \ref{figptw2}(c),(d) check the accuracy by displaying two types of errors as functions of the number of iterations:
\begin{itemize}
\item The residual error (in Euclidean norm)
\begin{eqnarray}
\label{fsec32}
RES_{\nu}=||L\varphi^{[\nu]}-N(\varphi^{[\nu]})||,\quad \nu=0,1,\ldots
\end{eqnarray}
\item
The discrepancy between the sequence of stabilizing factors  $s_{\nu}=s(\varphi^{[\nu]})$ and (in case of convergence) its limit one
\begin{eqnarray}
\label{fsec33}
SFE_{\nu}=|s_{\nu}-1|,\quad \nu=0,1,\ldots
\end{eqnarray}
\end{itemize} 
\subsection{Extended Boussinesq system}
The second example attempts the application of the method to systems of the form (\ref{ptw1}). Specifically we consider the so-called e-Boussinesq system
\begin{eqnarray}
\eta_{t}&=&-d_{1}W_{x}-d_{2}W_{xxx}-d_{4}(W\eta)_{x}+d_{5}(W\eta^{2})_{x},\label{ptw13a}\\
W_{t}&=&-\frac{1}{d_{1}}\eta_{x}-d_{3}W_{xxt}-\frac{d_{4}}{2}(W^{2})_{x}+d_{5}(W^{2}\eta)_{x},\label{ptw13b}
\end{eqnarray}
Equations (\ref{ptw13a}), (\ref{ptw13b}) appear as a model for the bidirectional propagation of an interfacial wave $\eta(x,t)$ between two fluid layers with $W(x,t)$ standing for the horizontal velocity of the flow. The constants are defined as
\begin{eqnarray*}
&&d_{1}=\frac{H}{r+H},\quad
d_{2}=\frac{H^{2}}{2(r+H)^{2}}(s+\frac{2}{3}(1+r H)),\nonumber\\
&&d_{3}=\frac{sd_{1}}{2},\quad
d_{4}=\frac{H^{2}-r}{(r+H)^{2}},\quad
d_{5}=\frac{r(1+H)^{2}}{(r+H)^{3}}, \label{ebou3}
\end{eqnarray*}
where $r$ and $H$ are, respectively, the (dimensionless) density and depth ratios while $s$ is related to the physical depth below and above the unperturbed interface and is taken within the range $-(1+rH)\leq s\leq -\frac{2}{3}(1+rH)$. The model is described in detail in \cite{NguyenD2008}, where the propagation and collision of solitary waves and fronts are also studied numerically. To our knowledge existence of periodic traveling wave solutions has not been analyzed (theoretically or numerically). Suggested here is some numerical evidence of such waves by using the method described in Section \ref{se2}.

Note first that periodic traveling waves $\eta=\eta(X), W=W(X), X=x-c_{s}t$ of (\ref{ptw13a}), (\ref{ptw13b}) of speed $c_{s}$ must satisfy the system for the profiles 
\begin{eqnarray}
\begin{pmatrix}
c_{s}&-(d_{1}+d_{2}\partial_{xx})\\
-\frac{1}{d_{1}}&c_{s}(1+d_{3}\partial_{xx})
\end{pmatrix}
\begin{pmatrix}
\eta\\W
\end{pmatrix}
=
\begin{pmatrix}
-W\eta(-d_{4}+d_{5}\eta)\\
-W^{2}(-\frac{d_{4}}{2}+d_{5}\eta)
\end{pmatrix}+\begin{pmatrix}
A_{1}\\A_{2}
\end{pmatrix},\label{ptw14}
\end{eqnarray}
for some integration constants $A_{1}, A_{2}$. The search for constant solutions $\eta=C_{1}, W=C_{2}$ of (\ref{ptw14}) leads to the algebraic system
\begin{eqnarray*}
\begin{pmatrix}
c_{s}&-d_{1}\\-\frac{1}{d_{1}}&c_{s}
\end{pmatrix}
\begin{pmatrix}
C_{1}\\C_{2}
\end{pmatrix}=\begin{pmatrix}
d_{4}C_{1}C_{2}-d_{5}C_{1}^{2}C_{2}\\\frac{d_{4}}{2}C_{2}^{2}-d_{5}C_{1}C_{2}^{2}
\end{pmatrix}+\begin{pmatrix}
A_{1}\\A_{2}
\end{pmatrix},
\end{eqnarray*}
in such a way that if $c_{s}^{2}\neq 1$ then
$$
C_{2}=\frac{c_{s}C_{1}-A_{1}}{d_{1}+d_{4}C_{1}-d_{5}C_{1}^{2}},
$$ while $C_{1}$ must be a root of the polynomial 
\begin{eqnarray*}
P(z)&=&-\frac{1}{d_{1}}z(d_{1}+d_{4}z-d_{5}z^{2})^{2}+c_{s}(c_{s}z-A_{1})(d_{1}+d_{4}z-d_{5}z^{2})\\
&&+(d_{5}z-\frac{d_{4}}{2})(c_{s}z-A_{1})^{2}-A_{2}(d_{1}+d_{4}z-d_{5}z^{2})^{2}.
\end{eqnarray*}
Finally, the system for the differences $\widetilde{\eta}=\eta-C_{1}, \widetilde{W}=W-C_{2}$ is of the form
\begin{eqnarray}
&&\left[\begin{pmatrix}
c_{s}&-(d_{1}+d_{2}\partial_{xx})\\
-\frac{1}{d_{1}}&+c_{s}(1+d_{3}\partial_{xx})
\end{pmatrix}
+
\begin{pmatrix}
(-d_{4}+2d_{5}C_{1})C_{2}&C_{1}(-d_{4}+d_{5}C_{1})\\
d_{5}C_{2}^{2}&2C_{2}\left(-\frac{d_{4}}{2}+d_{5}C_{1}\right)
\end{pmatrix}\right]
\begin{pmatrix}
\widetilde{\eta}\\\widetilde{W}
\end{pmatrix}\nonumber\\
&&=-
\begin{pmatrix}
d_{5}C_{2}\widetilde{\eta}^{2}
+\widetilde{W}\widetilde{\eta}(-d_{4}+d_{5}\widetilde{\eta})+d_{5}\widetilde{W}\widetilde{\eta}^{2}\\
\widetilde{W}^{2}(-\frac{d_{4}}{2}+d_{5}\widetilde{\eta})+2d_{5}C_{2}\widetilde{W}\widetilde{\eta}
+d_{5}\widetilde{W}^{2}\widetilde{\eta}
\end{pmatrix}.\label{ptw15}
\end{eqnarray}
Then, as mentioned in Section \ref{se2}, system (\ref{ptw15}) is iteratively solved for $\widetilde{\eta}, \widetilde{W}$ and approximate periodic traveling wave profiles for the values $A_{1}, A_{2}$ are obtained from the formulas $\eta=\widetilde{\eta}+C_{1}, W=\widetilde{W}+C_{2}$.

In \cite{NguyenD2008} a maximum wave velocity
$$
v_{max}=1+\frac{(H^{2}-r)^{2}}{8rH(1+H)^{2}},
$$ is derived and solitary waves are characterized by wave velocities larger than one. We illustrate the performance of the method by generating periodic traveling waves with a speed close to $v_{max}$ (specifically $c_{s}=v_{max}-10^{-4}$ in the experiments below) and two different values for $r, H$ and the constants $A_{1}, A_{2}$. (In both experiments $s=-(1+rH)$ has been fixed.) 
\begin{figure}[htbp]
\centering 
\subfigure[]{
\includegraphics[width=6.6cm]{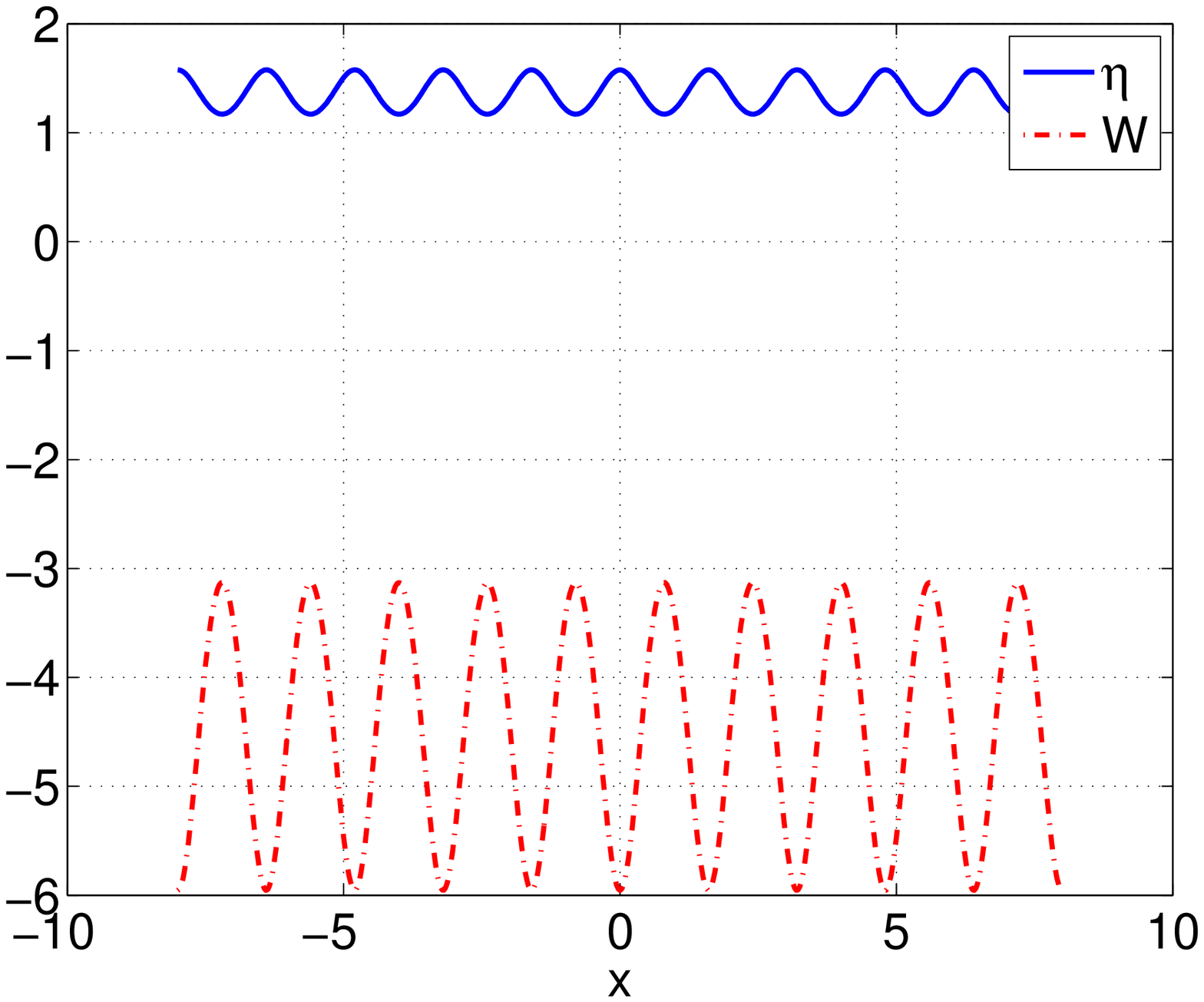} }
\subfigure[]{
\includegraphics[width=6.6cm,height=5.6cm]{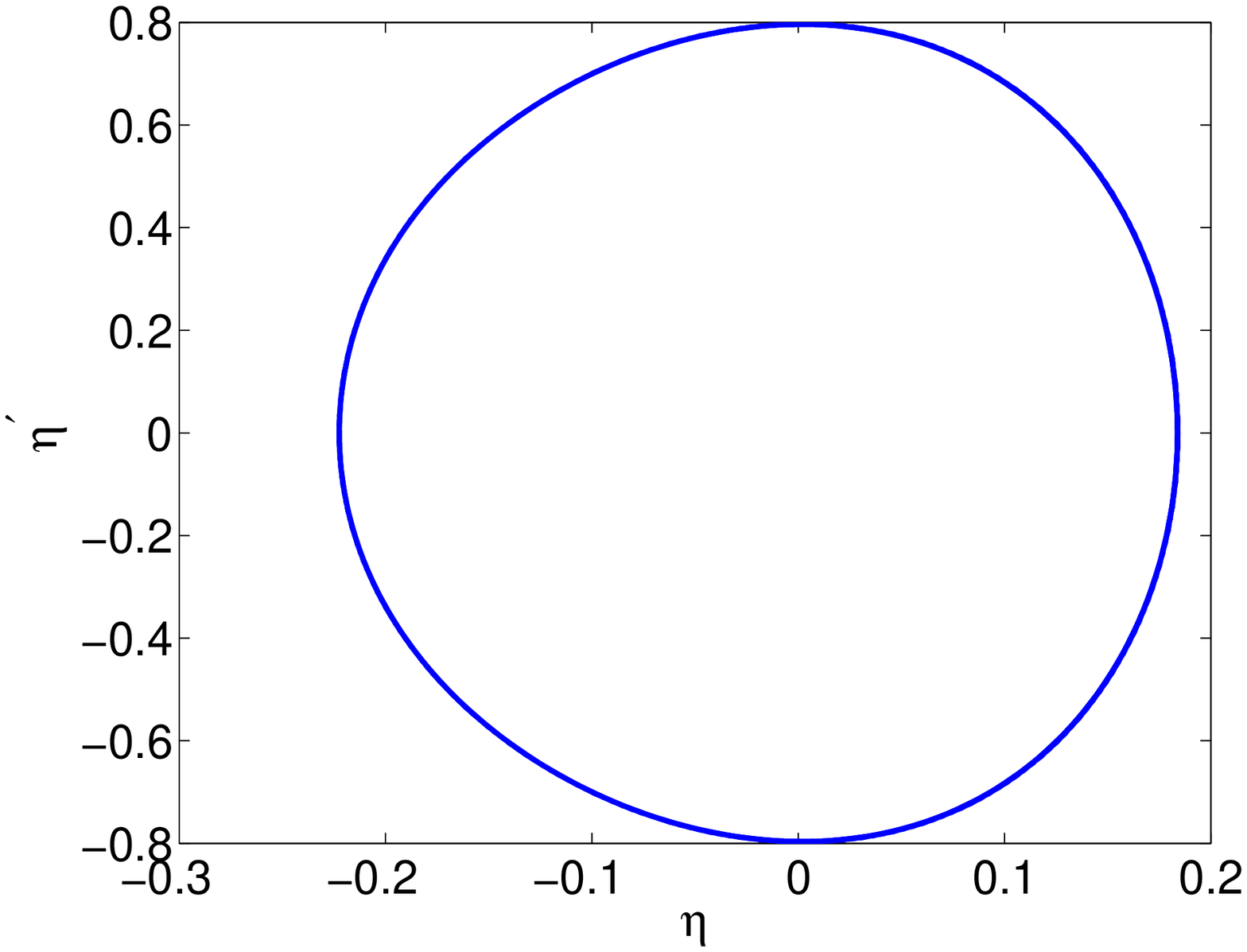} }
\subfigure[]{
\includegraphics[width=6.6cm]{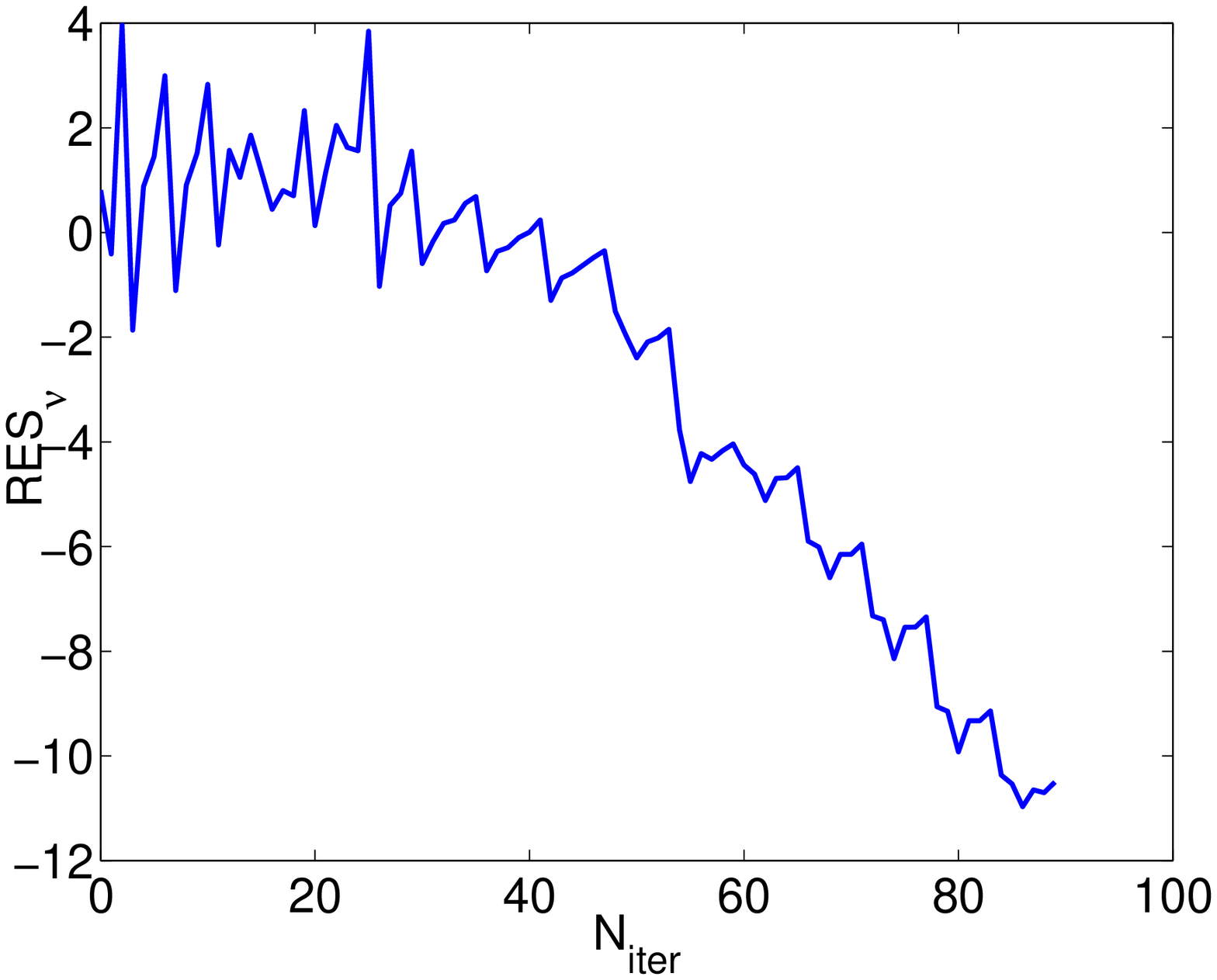} }
\subfigure[]{
\includegraphics[width=6.6cm,height=5.4cm]{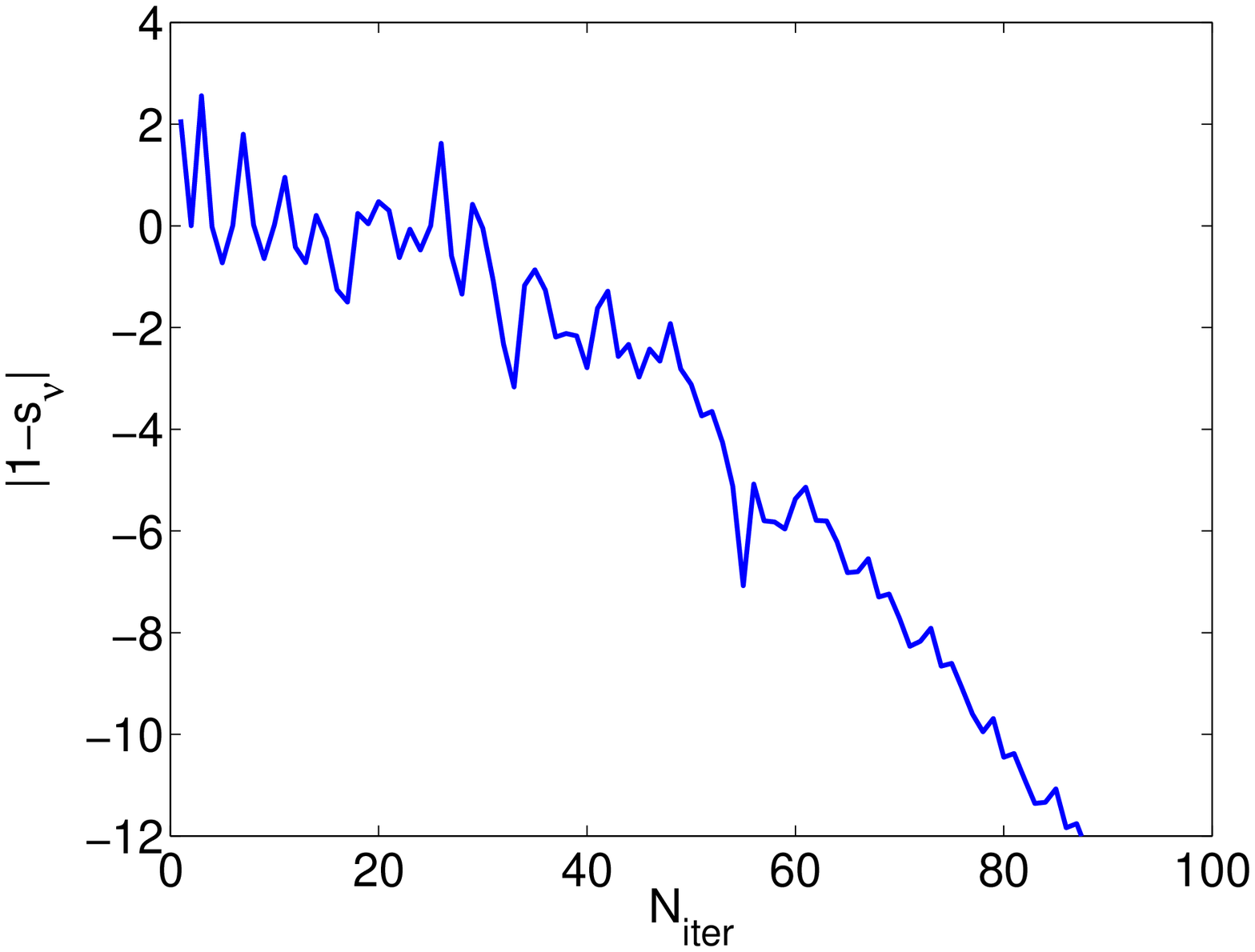} }
\caption{Periodic traveling wave generation of  (\ref{ptw14}) with $c_{s}=v_{max}-10^{-4}$. $A_{1}=-1, A_{2}=-2, r=0.8, H=0.95, s=-1-rH$. $v_{max}=1.000454$. (a) Approximate profiles; (b) $\eta$ phase portrait; (c) Residual error (\ref{fsec32}) vs number of iterations; (d) Discrepancy (\ref{fsec32}) vs number of iterations (semi-logarithmic scale).} \label{figptw3}
\end{figure}
The first case, Figure \ref{figptw3}, corresponds to $r=0.8, H=0.95, A_{1}=-1, A_{2}=-2$. Figure \ref{figptw3}(a) shows the approximate $\eta, W$ profiles (also obtained by using additional acceleration with the MPE method) and the periodic character is confirmed by the phase portrait of $\eta$ in Figure \ref{figptw3}(b). The accuracy of the computation is again evaluated in Figures \ref{figptw3}(c) and \ref{figptw3}(d) with the behaviour of the errors (\ref{fsec32}) and (\ref{fsec33}) with respect to the number of iterations. 

\begin{figure}[htbp]
\centering \subfigure[]{
\includegraphics[width=6.6cm]{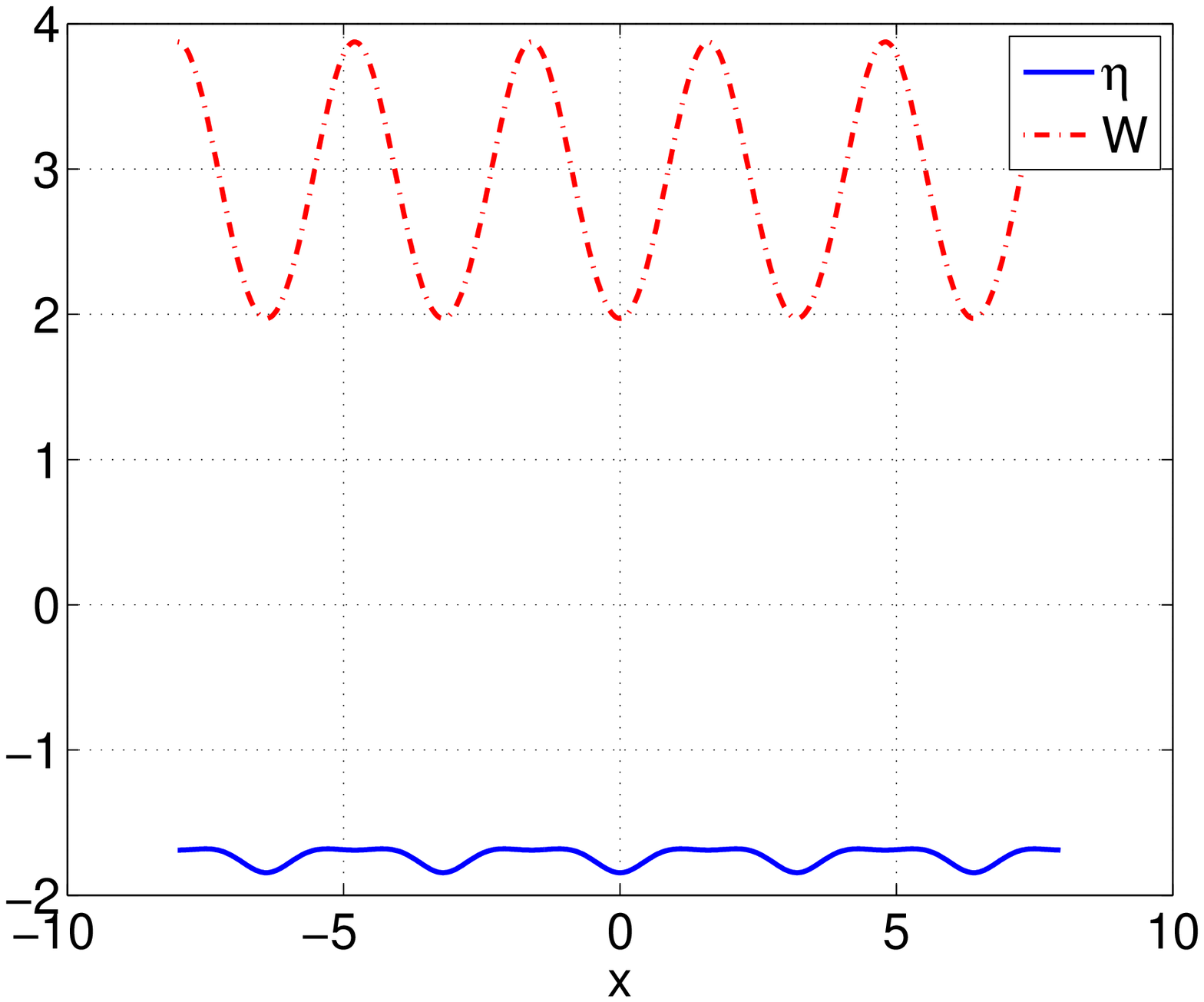} }
\subfigure[]{
\includegraphics[width=6.6cm,height=5.6cm]{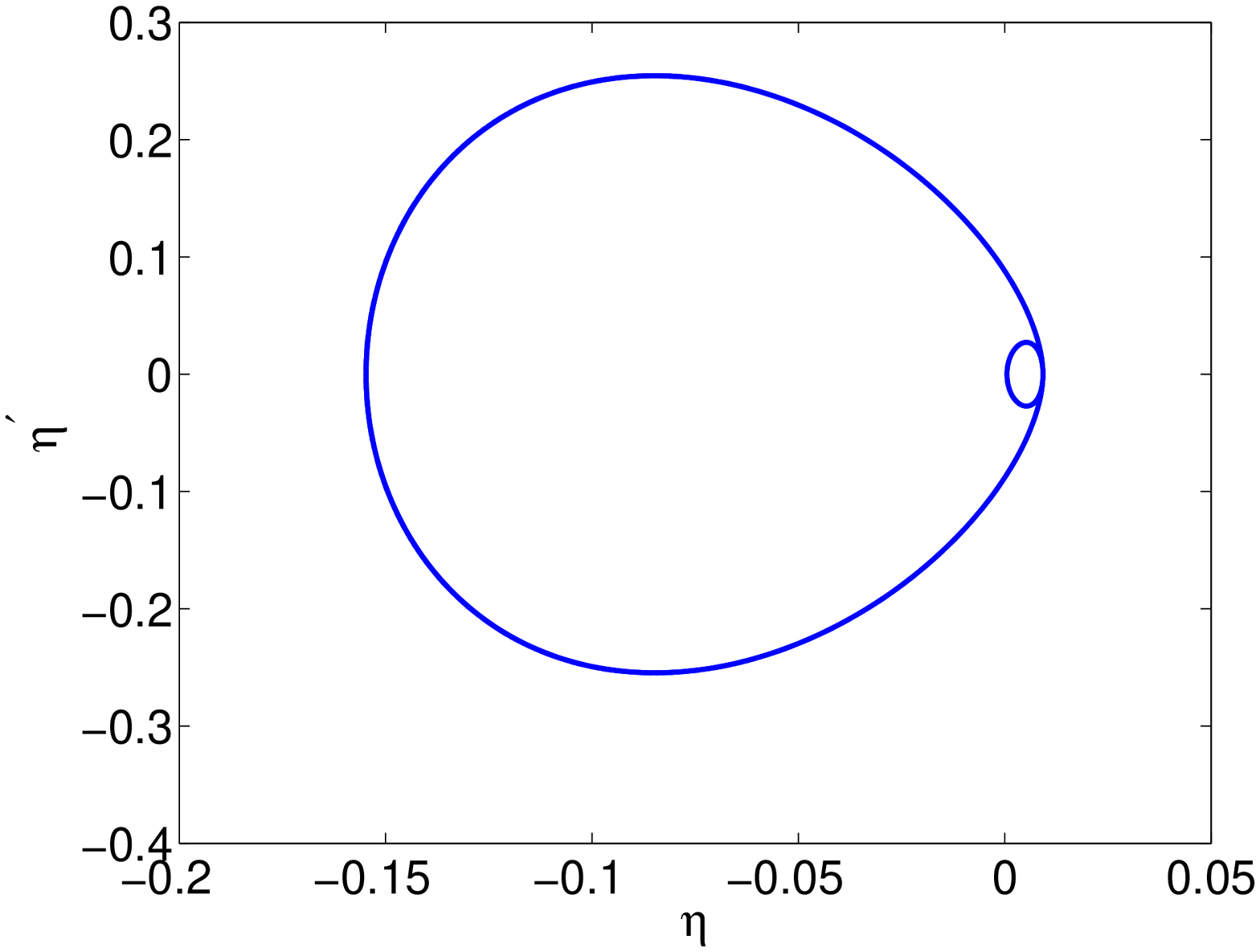} }
 \subfigure[]{
\includegraphics[width=6.6cm]{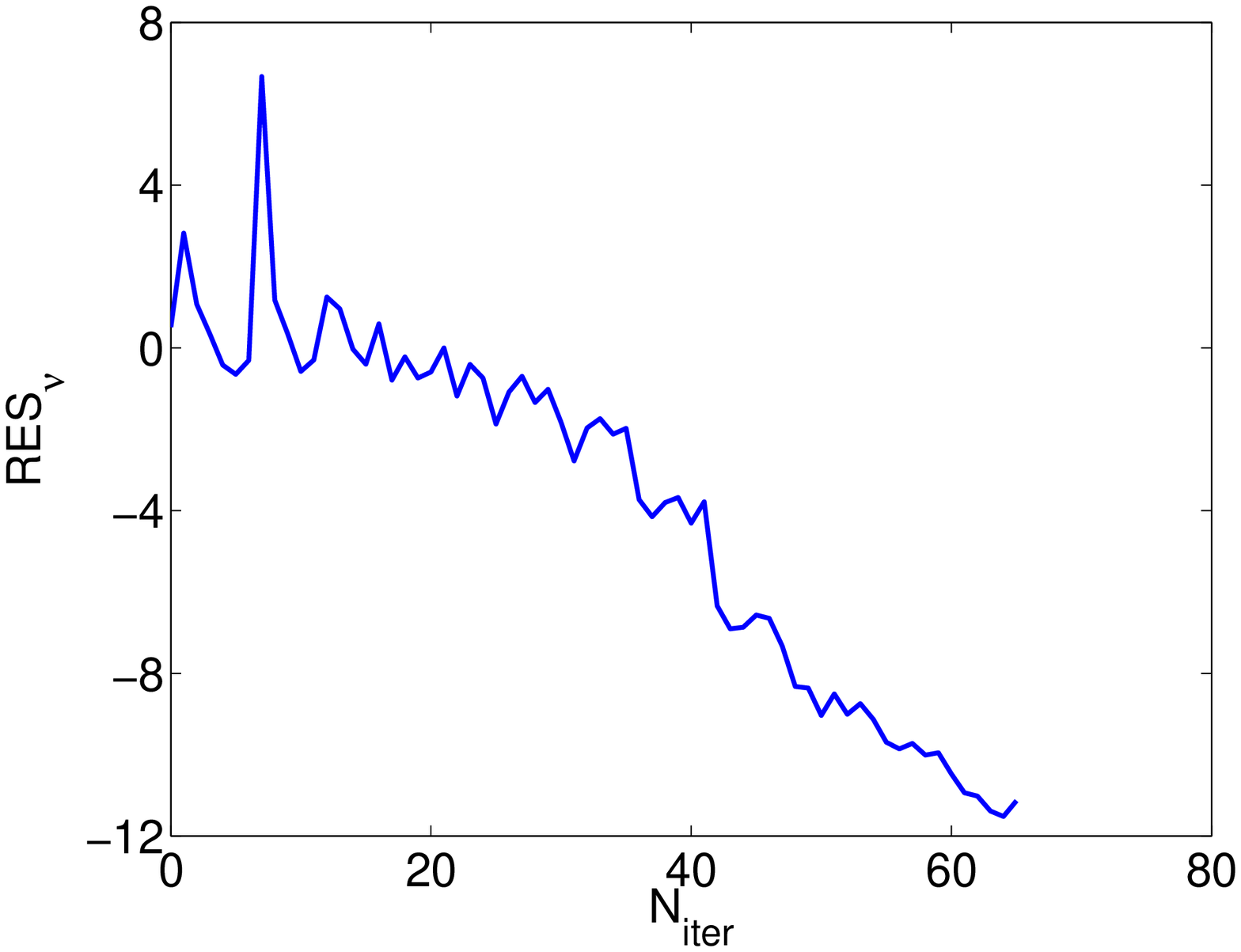} }
\subfigure[]{
\includegraphics[width=6.6cm]{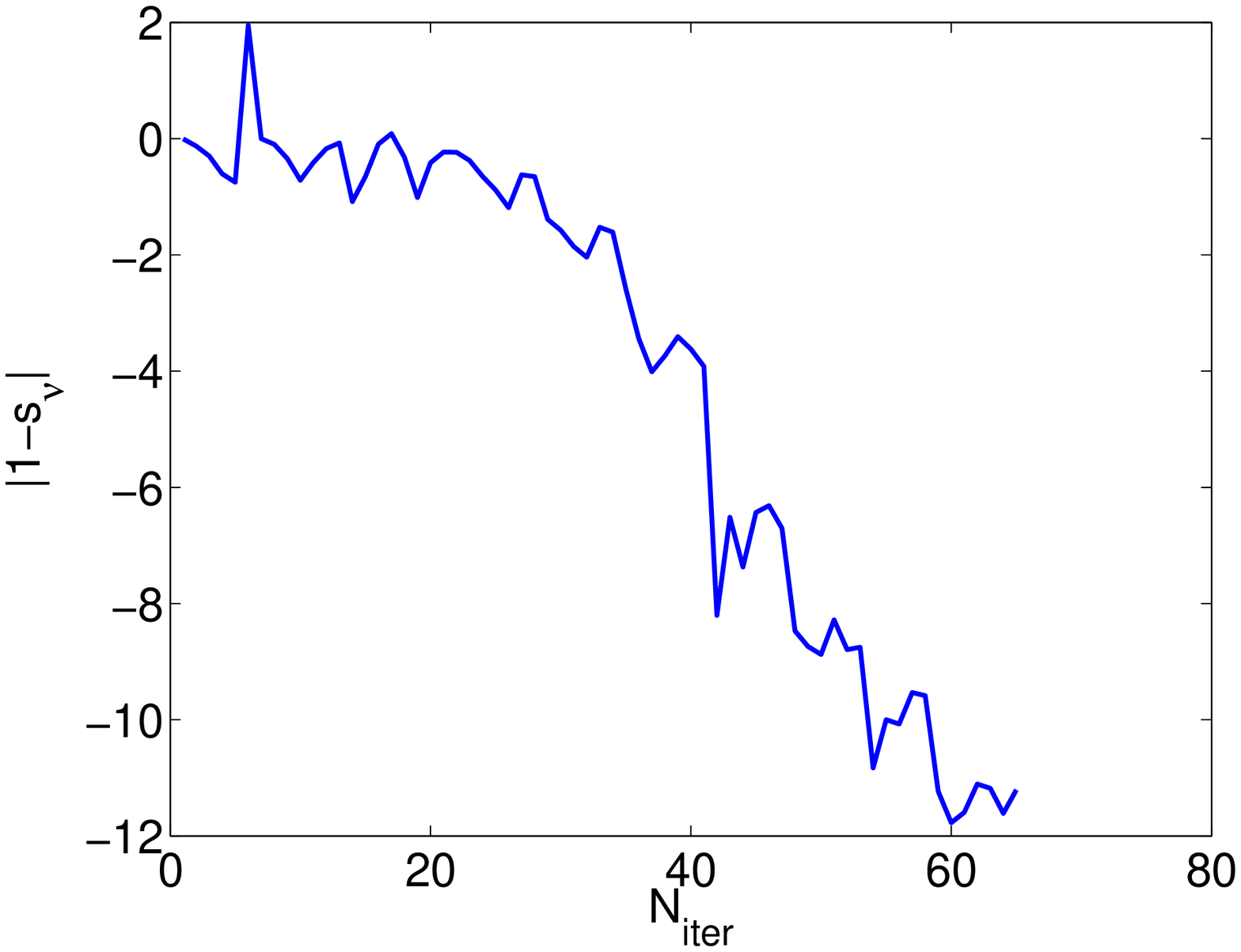} }
\caption{Periodic traveling wave generation of  (\ref{ptw14}). Approximate profiles and phase portraits with $c_{s}=v_{max}-10^{-4}$. $A_{1}=1, A_{2}=1, r=0.8, H=1.8, s=-1-rH$. $v_{max}=1.065919$. (a) Approximate profiles; (b) $\eta$ phase portrait; (c) Residual error (\ref{fsec32}) vs number of iterations; (d) Discrepancy (\ref{fsec32}) vs number of iterations (semi-logarithmic scale). } \label{figptw4}
\end{figure}

The same information is provided in Figure \ref{figptw4} for the case $r=0.8, H=1.8, A_{1}=1, A_{2}=1$. Note here that the periodic form of the approximate $\eta$ profile is different and closer to the \lq table-top\rq\ form of the solitary waves shown in \cite{NguyenD2008}. This suggests the connection between these periodic profiles and some associated solitary wave in the limit of large wavelength, \cite{Bona1981}.
\section{Concluding remarks}
\label{se4}
In this paper an iterative technique to generate periodic traveling wave solutions of some nonlinear dispersive wave equations is proposed. The method is based on a suitable change of variable to avoid the presence of some constants of integration in the problem and the application of fixed point iteration of Petviashvili type, \cite{AlvarezD2014b}. This strategy overcomes two difficulties not totally solved in some previous algorithms, \cite{AlvarezD2011}, for this purpose, namely:
\begin{itemize}
\item The presence of constants that prevents the use of Petviashvili type methods, which have been shown to be an alternative to compute solitary wave solutions for the models under study (see \cite{AlvarezD2014a} and references therein).
\item The presence of nonlinearities of polynomial type with homogeneous terms of different degree, for which the construction of the Petviahsvili type methods is not suitable.
\end{itemize}
The scheme has been illustrated by computing periodic traveling waves in two problems of additional interest: the generalized fractional KdV equation and the extended Boussinesq system. In the first case, existence and stability of periodic traveling waves have been studied theoretically but, to the authors' knowledge, no approximate profiles have been computed. The second example attempts to extend the application of the technique to some nonlinear dispersive systems and suggests the formation of periodic traveling waves in the model, for which (again to the authors' knowledge) there are no existence results.
\section*{Acknowledgements}
This research has been supported by  project MTM2014-54710-P.

\end{document}